\documentclass{llncs}

\usepackage{url}
\usepackage{cite}
\usepackage{graphicx}
\usepackage{framed}
\usepackage{amsfonts}
\usepackage{amsmath}

\begin{document}

\title{Another (wrong) construction of $\pi$}
\author{Zolt\'an Kov\'acs}
\institute{
The Private University College of Education of the Diocese of Linz\\
Salesianumweg 3, A-4020 Linz, Austria\\
\email{zoltan@geogebra.org}
}

\maketitle              

\begin{abstract}
A simple way is shown to construct the length $\pi$ from the unit length
with 4 digits accuracy.

\end{abstract}

\section{Introduction}
It is well-known that accurate construction of the ratio between the perimeter and
the diameter of a circle is theoretically impossible \cite{Wantzel1837,Lindemann1882}.
Before Lindemann's result (and even thereafter), however, several attempts
were recorded to geometrically construct the number $\pi$ with various means,
in most cases by using a compass and a straightedge. One of the most successful
attempts is Kocha\'nski's work \cite{Kochanski} that produces an approximation
of $\pi$ by 4 digits, $\sqrt{40/3-2\sqrt3}\approx3.14153333\ldots$

Kocha\'nski's construction is relatively easy, it requires just a little amount
of steps, and can be discussed in the school curriculum as well.
In this note the same approximation is given, by using---at least geometrically---an
even simpler approach.

\section{The construction}

\begin{figure}
\centering
\includegraphics[width=0.6\textwidth]{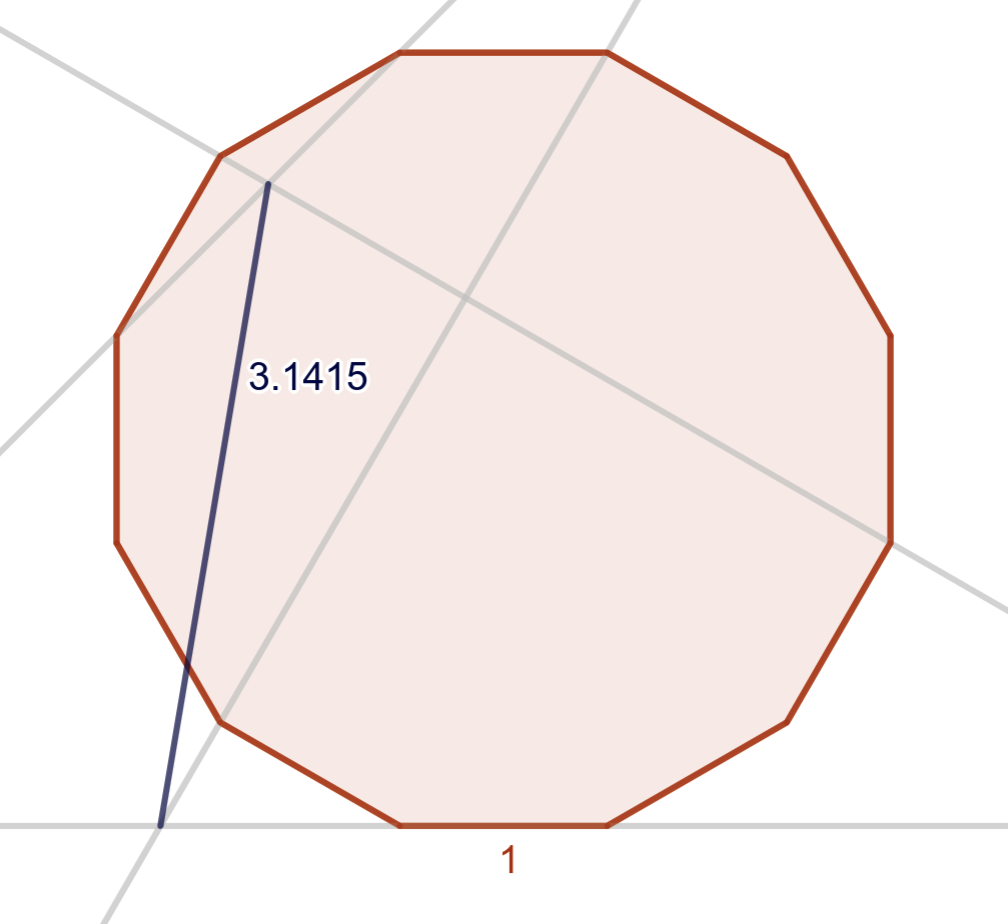}
\caption{A new way to construct $\pi$ approximately}
\label{fig1}
\end{figure}

The proposed way to construct $\pi$ is shown in Fig.~\ref{fig1} (see also \cite{pi-12gon}).

\begin{figure}
\centering
\includegraphics[width=0.6\textwidth]{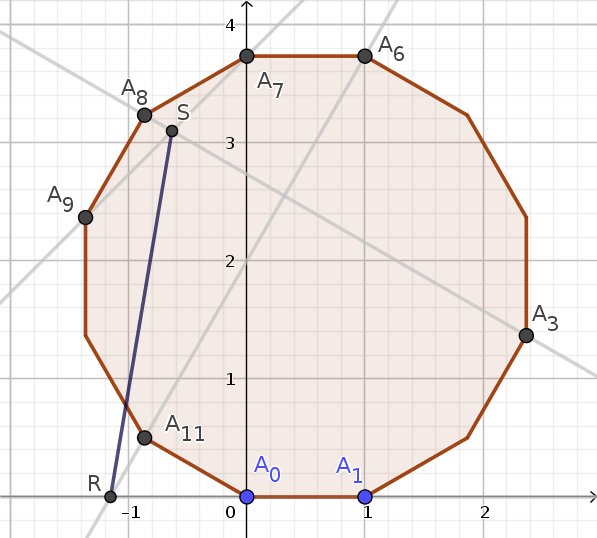}
\caption{Explanation for the proof}
\label{fig2}
\end{figure}

A proof that $|RS|=\sqrt{40/3-2\sqrt3}$ is as follows. We assume that $A_0=(0,0)$ and
$A_1=(1,0)$ (see Fig.~\ref{fig2}). By using that $360^{\rm o}/12=30^{\rm o}$,
$\cos(30^{\rm o})=\sqrt3/2$ and $\sin(30^{\rm o})=1/2$, the exact coordinates of the appearing vertices in the construction are
$A_3=(3/2+\sqrt3/2,1/2+\sqrt3/2)$, $A_6=(1,2+\sqrt3)$, $A_7=(0,2+\sqrt3)$, $A_8=(-\sqrt3/2,3/2+\sqrt3)$, $A_9=(-1/2-\sqrt3/2,3/2+\sqrt3/2)$,
$A_{11}=(-\sqrt3/2,1/2)$.

To find the coordinates of $R$ we can compute the equation of the line $A_{11}A_6$.
By substituting the coordinates of $A_{11}$ and $A_6$ it can be verified that
the equation is $$(3/2+\sqrt3)x-(1+\sqrt3/2)y=-2-\sqrt3,$$
and solving this for $y=0$ we obtain the exact coordinates $R=(-2/\sqrt3,0)$.
(Alternatively, it can be shown that $|A_1R|=1+2/\sqrt3$, because it is
the shorter cathetus of the triangle $RA_1A_6$ which is a half of an equilateral triangle---this
holds because the angle $A_1A_6A_{11}$ is an inscribed angle of the circumcircle of
the regular 12-gon, and it must be $60^{\rm o}/2=30^{\rm o}$.)

Now, $A_3A_8\perp A_6A_{11}$, so we are searching for the equation of line $A_3A_8$ in form
$$(1+\sqrt3/2)x+(3/2+\sqrt3)y=c.$$ After substituting the coordinates of $A_3$ in this,
we obtain that $c=(9+5\sqrt3)/2$. Because of symmetry, it is clear that the equation of line
$A_7A_9$ is of form $$y=x+d.$$ By using the coordinates of $A_7$, we immediately obtain
that $d=2+\sqrt3$.

To find the coordinates of $S$ we solve the equation system
\begin{align}
(1+\sqrt3/2)x+(3/2+\sqrt3)y=&(9+5\sqrt3)/2,\\
y=&x+2+\sqrt3
\end{align}
now, which produces the coordinates $x=(\sqrt3-3)/2$, $y=(3\sqrt3+1)/2$.

Finally we compute the length of $RS$:
$$|RS|=\sqrt{\left((\sqrt3-3)/2+(2/\sqrt3)\right)^2+\left((3\sqrt3+1)/2\right)^2}=\sqrt{40/3-2\sqrt3}.$$

\section{Remarks}

This result has been found by the software tool \cite{RegularNGons} in an automated way by considering all possible
configurations of distances between intersections of diagonals in regular polygons.
It seems very likely that the construction described above is the simplest and most
accurate one among the considered cases. Another construction which is based on a regular \textit{star}-12-gon
(see \cite{pi-star12gon}) can produce the same approximation.

Regular polygons with less
sides have already been completely studied for the constructible cases for $n<12$ with
less accurate results (see \cite{pi-12gon} for details). Checking cases $n=15,16,17,20$
(all are constructible) is an on-going project.

While the software tool \cite{RegularNGons} gave a machine assisted proof by using Gr\"obner bases
and elimination (see \cite{Cox_2007,Recio1999} for more details), the proof above was compiled by
the author manually.

\section*{Acknowledgments}
The author was partially supported by a grant MTM2017-88796-P from the
Spanish MINECO (Ministerio de Economia y Competitividad) and the ERDF
(European Regional Development Fund).

\bibliography{kovzol,external}

\end{document}